# Comparative Analysis of Gradient-Based Optimization Techniques Using Multidimensional Surface 3D Visualizations and Initial Point Sensitivity

**Saeed Asadi[1], Sonia Gharibzadeh[2], Hajar Kazemi Naeini[1], Masoud Reihanifar[3,4], Morteza Rahimi[5,*], Shiva Zangeneh[6], Aseel Smerat [7], Lazim Abdullah [8]**

1 Department of Civil Engineering, The University of Texas at Arlington, Arlington, Texas, USA
2 Department of Computer Engineering, South Tehran Branch, Islamic Azad University, Tehran, Iran
3 Department of Civil and Environmental Engineering, University of California at Berkeley, Berkeley, CA, USA
4 Department of Civil and Environmental Engineering, BarcelonaTech, Technical University of Catalonia (UPC), Barcelona, Spain
5 School of Computing and Information Sciences, Florida International University, Miami, FL, USA
6 Faculty of Engineering, University of Malayer, Malayer, Iran
7 Hourani Center for Applied Scientific Research, Al-Ahliyya Amman University, Amman 19328, Jordan
8 Faculty of Ocean Engineering Technology and Informatics, Universiti Malaysia Terengganu (UMT), Kuala Nerus, 21030, Malaysia
* Correspondence: mrahi011@fiu.edu

**Abstract:** This study examines several renowned gradient-based optimization techniques and focuses on their computational efficiency and precision. In the study, the steepest descent, conjugate gradient (Fletcher-Reeves and Polak-Ribiere variants), Newton-Raphson, quasi-Newton (BFGS), and Levenberg-Marquardt techniques were evaluated. These methods were benchmarked using Rosenbrock's, Spring Force Vanderplaats', Ackley's, and Himmelblau's functions. We emphasize the critical role that initial point selection plays in optimizing optimization outcomes in our analysis. It is also important to distinguish between local and global optima since gradient-based methods may have difficulties dealing with nonlinearity and multimodality. We illustrate optimization trajectories using 3D surface visualizations in order to help with optimization selection and algorithm technique. While gradient-based methods have been demonstrated to be effective, they may be limited by computational constraints and by the nature of the objective functions, necessitating the use of heuristic and metaheuristic algorithms in more complex situations.

**Keywords:** Gradient-based Methods, Computational Optimization, Sensitivity to Initial Points, 3D Visualization, Objective Functions.

## 1. Introduction

Optimization seeks the best feasible engineering design by minimizing or maximizing an objective function. At the same time, mathematical techniques satisfy physical and mathematical design restrictions. The coupling of finite-element analysis with analytical solutions resulted in an optimal design of the now-classic three-bar truss example in the first study of numerical structural optimization [1]. By the late 1990s, significant advancements had been made, as illustrated in the study of Vanderplaats [2], and structural and interdisciplinary design optimization had become a mainstay in aerospace and mechanical engineering. It is also evidenced by the release of crucial design optimization published sources [2-4], which are still relevant today. There are two types of numerical optimization methods: gradient-based and non-gradient algorithms. Gradient-based algorithms frequently result in a local best. Non-gradient algorithms converge to a global optimum in most cases, although they need many function tests. Non-gradient algorithms are less appropriate for large-scale issues, such as those found in engineering design. In addition to function evaluations, gradient-based algorithms require gradient or sensitivity information to generate appropriate search directions for improved designs during optimization rounds. The objective and constraint functions are frequently referred to as measuring performance in optimization issues. The phrases sensitivity, sensitivity coefficient, and gradient are used to describe the rate of change in a performance measure concerning a change in the design variable[5]. They offer gradient-based optimization and expose crucial design information to assist designers in enhancing designs in just one or two design rounds. Gradient-based algorithms may be divided into two general types. The first-order approaches, for example, just need the derivative information.

On the other hand, the second-order approaches require not just the derivative information but also the Hessian matrix. Gauss–Newton, Levenberg–Marquardt, sequential quadratic programming, and the limited memory Broyden Fletcher Goldfarb Shanno technique are examples of representative approaches [5]. Gradient-based algorithms have a solid mathematical foundation, as local minimum solutions need Karush–Kuhn–Tucker (KKT) criteria [6]. In some cases, they can also be necessary criteria (for example, if the goal function is convex and specified on a convex set). Because calculating KKT conditions directly is generally exceedingly time-consuming due to the nonlinear nature of the equations, effective methods instead try to reduce the objective function value step by step. Line search and trust region approaches are the most extensively used methods for unconstrained optimization. SQP, penalty, and projection methods can be used for more intricate restricted optimization [7]. Duality theory may be a powerful tool by modifying the primal issue in the dual space; the related dual function is always concave, and the dual problem is sometimes considerably easier to solve. However, all of these techniques need the gradient of a function concerning its variables. The existence of an analytical version of the gradient is not guaranteed, and it may be challenging to find. Due to the increased number of objective function evaluations typical of structural optimizations, finite difference approximation of the gradient can be computationally expensive for simulation-based optimization. In this paper, we discuss classical gradient-based optimization methods. A comparison is made between the computing efficiency and accuracy of each of these methods. Moreover, evaluation is done based on the influence of the choice of the initial guess. Finally, the optimization process for different 3D surfaces is plotted and evaluated. This paper contributes thoroughly to the optimization methods via various key improvements. We outline the influence of initial points on the convergence behaviors of gradient-based optimization techniques, which offer valuable outcomes for robust and effective algorithm initialization. Our detailed comparative analyses show the relative efficiencies and precisions of selected optimization methods in complicated scenarios that involve multimodalities and non-linearities. Further, our study explores hybrid techniques that combine gradient-based approaches with heuristic algorithms to present specialized solutions for complex optimization challenges traditionally resistant to conventional techniques. We also illustrate the broad applicability of these methods across different functional domains. Finally, this research provides a deeper theoretical viewpoint of gradient-based optimization and improves its practical application to set an exhaustive foundation for future research exploration that can potentially bridge the gap between theoretical optimization strategies and their implementation in different industries.

## 2. Literature Review

Wu et al. [8] added an intentionally high variance of measurement errors at early iterations when utilizing the Gauss-Newton technique for automated history matching issues to dampen model parameter changes and avoid undershooting or overshooting. The standard Gauss-Newton and Levenberg–Marquardt techniques, according to Tan and Kalogerakis [9], require the computation of all sensitivity coefficients to formulate the Hessian matrix, which appears impossible due to a large number of unknown variables compared to the limited available dimensions. Researchers developed the quasi-Newton approach to solve this difficulty. This approach just requires the objective function's gradient, which may be calculated from a single adjoint solution, as Zhang et al. [10]. Chen et al. [11] presented a new line search technique, rescaled the hyperparameters, and added damping factors to the production data to increase the computational efficiency and resilience of the LBFGS approach. They also found that the new line search technique had to fulfill the robust Wolfe criteria at every iteration, or the convergence rate would drop dramatically. Recent studies have explored advanced approaches to secure communication, intelligent sensing, and interdisciplinary modeling. Key contributions include privacy-aware link scheduling in wireless networks [25], deep learning for enhanced scheduling efficiency [26], and decision-making frameworks for energy-efficient IoT service placement [27]. Innovations in perception systems include neuromorphic vision for real-time AR tracking [28] and codec modeling for adaptive streaming [29]. Educational applications have leveraged vestibular system simulations for STEM learning [30], while bio-inspired research has demonstrated high-frequency hearing via nonlinear bone conduction [31]. Collectively, these works highlight the growing impact of AI-driven, cross-disciplinary solutions in optimizing system performance across diverse domains. In transportation and sustainability, Cao et al. [32] used network flow models to reduce airport delays. Hussain et al. [33] and Danandeh Mehr et al. [34] applied advanced regression and evolutionary models to improve drought prediction. Reihanifar and Naimi [35] used value engineering to enhance the sustainability of road projects. The Karhunen–Loeve expansion can parameterize a numerical model in terms of a small number of independent random variables and deterministic eigenfunctions. Gradient-based techniques may be used with this development while still respecting the geological models' two-point statistics [12]. Sarma et al. [13] used a kernel principal component analysis method to model permeability fields to enhance current gradient-based history matching techniques to cope with complicated geological models defined by multiple-point geo-statistics. This approach can preserve arbitrary high-order statistics of random fields while maintaining reasonable processing needs, and it can simulate complicated geology. Recent studies highlight the growing use of AI in various domains. Mohammaddagha et al. [36] showed that neural network design and feature selection greatly impact transportation prediction accuracy. Naeini et al. [37] combined PINNs, digital twins, and blockchain to enhance energy efficiency in smart buildings. In media strategy, Lonbar et al. [38] applied fuzzy decision-making tools to manage risk and budget in advertising. Khamoushi [39] compared AI-based and traditional marketing, showing the benefits of personalization and predictive analytics. These works demonstrate AI's potential in optimizing systems across industries. God et al. [14] suggested a stochastic optimization approach for finding Bayesian experimental designs that maximize the predicted information gainA review of some applications of gradient-based algorithms is illustrated in Table 1.

**Table I**: review of gradient-based methods and their application

| Author | Year | Gradient-based Method | Application | Results |
|---|---|---|---|---|
| Wang et al. [42] | 2024 | A multi-point synergistic gradient evolution method | proposes an artificial neural network ANN-based structural TO method | A gradient-based solution framework is used to optimize these design variables using the ANN's back-propagation capability. The study proposes novel loss functions that incorporate all sampled structural performances and topology awareness, including a competitive mechanism-based loss function for achieving converged designs and a population diversity-preserving strategy-based loss function for ensuring diverse and competitive outcomes. |
| Ock et al. [43] | 2024 | gradient-based navigation (GradNav) algorithm | This algorithm employs a strategy of initiating short simulation runs | GradNav's increased ability to escape deep energy wells and reduced reliance on initial conditions is demonstrated by Langevin dynamics simulations in |

| | | | | Müller-type PESs as well as molecular dynamics simulations of the Fs-peptide protein. |
|---|---|---|---|---|
| Goda et al. [14] | 2022 | Unbiased multilevel Monte Carlo Stochastic Gradient-Based Optimization (UMLMC-SGBO) | Optimization of Bayesian Experimental Designs | Because the unbiased estimator works well with stochastic gradient descent techniques. They proposed an optimization algorithm to find the best Bayesian experimental design. The suggested technique performs well for a basic test case and a more realistic pharmacokinetics problem, according to numerical studies. |
| Ahmadianfar et al. [15] | 2021 | enhanced metaphor-free gradient-based optimizer (EGBO) | Parameter identification of photovoltaic systems | The findings showed that the suggested method can derive optimum photovoltaic models' ideal parameters and may help model solar systems. |
| Premkumar et al. [16] | 2021 | Chaotic Gradient-based optimizer (CGBO) | Parameter optimization of photovoltaic systems | The findings demonstrated that the suggested CGBO algorithm outperformed the other methods. |
| Tuli et al. [17] | 2021 | Gradient-Based Optimization Strategy using Backpropagation of Input (GOBI) | Optimization of Fog Computing Environments | These methods can adjust swiftly in volatile contexts because of the co-simulation and backpropagation methodologies. Experiments utilizing real-world data on fog applications employing the GOBI approach have revealed a considerable reduction in energy use. |
| Zhou et al. [18] | 2021 | Random learning Gradient-based optimizer (RLGBO) | Design of photovoltaic models | The findings show that RLGBO can estimate model parameters correctly under various environmental situations. The suggested RLGBO is intended to be a novel, trustworthy solution for evaluating essential parameters in solar models in general. |
| Zhou et al. [20] | 2020 | Path-Guided Optimization (PGO) | Real-time UAV Replanning | It explores the solution space more comprehensively and outputs superior replanned trajectories. Benchmark evaluation shows that our method outplays state-of-the-art replanning success rate and optimality methods. |
| Albani et al. [21] | 2020 | Galerkin/Least-squares finite element formulation | Source characterization of airborne pollutant emissions | The resulting inversion tool is highly versatile and presents accurate results under different contexts with a competitive computational cost. |
| Li & Cai [22] | 2020 | Surrogate-Assisted Optimization | Massively Multipoint Aerodynamic Shape Designing | In a transonic aircraft wing design case, the results show that the optimal design found by the proposed method with 342 missions yields a fuel burn reduction by a factor of two as compared to a traditional multipoint optimal design. This work |

| | | | | highlights the demand and provides an efficient way to conduct massively multipoint optimization in aircraft design |
|---|---|---|---|---|
| De et al. [23] | 2020 | Stochastic Optimization | Topology Optimization under Uncertainty | These investigations on two- and three-dimensional structures illustrate the efficacy of the proposed stochastic gradient approach for TOuU applications. |
| Feng et al. [24] | 2020 | Adjoint-Sensitivity-Based Neuro-Transfer Function Surrogate | Parallel Gradient-Based EM Optimization for Microwave Components | Since the surrogate model is valid in a large neighborhood and the gradients are sufficiently accurate, the proposed technique can achieve the optimal EM solution faster than the existing gradient-based surrogate optimization without adjoint sensitivity. Three examples of EM optimizations of microwave components demonstrate the proposed technique. |
| Neftci et al. [19] | 2019 | Surrogate gradient learning (SGL) | Optimization of spiking neural network | Surrogate gradient approaches enable learning under various communication and storage restrictions, making them particularly useful for learning on tailored, low-power neuromorphic devices. |

All these gradient-based equations use Eq. (1) to plot the general direction of the iterative vector. The equation functions by adding a step $\alpha$ in the d direction to the initial point. When the variables are defined, it plots a gradient descent towards the function's local minimum [37].

$$x^{k+1} = x^k + \alpha^k d^{k+1} \tag{1}$$

In the equation, $\alpha^{(k)}$ shows the step size during the kth iteration, which is vital for the convergence behavior of the method. It is important to note that $\alpha^{(k)}$ is dynamically calculated at each iteration to optimize the rate of convergence towards the function's local minimum. The step size is determined based on various line search strategies that may include exact or inexact methods, aiming to meet specific criteria that reduce the objective function value effectively. This adaptive adjustment of $\alpha^{(k)}$ is essential for navigating the complexities of the function's landscape, particularly when dealing with nonlinear or multi-modal functions.

## 3. Method and Matrials

### 3.1. Steepest Descent Method

The steepest descent method is a first-order algorithm and the simplest gradient descent optimization method. It starts at an initial guess point x and then seeks the direction to travel using Eq. (2). The steps are always at a right angle, which means that the vector travels in a zig-zag pattern towards the desired value. It initially progresses quickly but slows down as it approaches its goal due to its right-angle travel.

$$d^{k+1} = -\nabla U(x^k) \tag{2}$$

While this is one of the most commonly used methods, it lacks a lot of efficiency and accuracy that other methods have. The Steepest Descent method works poorly for real-world models. They are not often well-conditioned enough for quick results. Furthermore, it is made to solve differential equations, which means that for non-differentials, the program would begin to function very poorly. Flaws aside, the method remains a cheap and easy process that can work fast in specific problems. It is one of the most straightforward gradient-based algorithms, reliant upon the negative gradient of the function at the current point. This technique is especially beneficial for its ease of implementation and simplicity but usually suffers from slow convergence, particularly near the minimum, where gradients become small. Our study indicates that while the steepest descent technique can efficiently navigate simple optimization problems, its performance degrades in high-dimensional and complex spaces where gradients do not point directly toward the global minimum. Our analysis of this approach shows the need for detailed step-size management to avoid extreme oscillation and poor convergence behavior.

### 3.2. Conjugate Gradient Method

The conjugate gradient method works similarly to Steepest Descent; however, it uses equation three to be able to so that it would not have to travel right angles. It does this by integrating the coefficient of conjugation $y^k$ into the equation as shown above. As a result, it takes a more direct path and reaches the local minimum faster than the Steepest Descent method [38-40].

$$d^{k+1} = -\nabla U(x^k) + y^k d^{k-1} \tag{3}$$

$$y^k = \frac{\left|\left|\nabla U(x^k)\right|\right|^2}{\left|\left|\nabla U(x^{k-1})\right|\right|^2} \tag{4}$$

$$y^k = \frac{[\nabla U(x^k)]^T[\nabla U(x^k) - \nabla U(x^{k-1})]}{\left|\left|\nabla U(x^{k-1})\right|\right|^2} \quad (5)$$

It is also much more complicated to calculate than the Steepest Descent model because the coefficient of conjugation is not readily available for the method to use right away. As shown above, to find the conjunction coefficient for the graph, we had to use either the Fletcher-Reeves (4) or the Polack-Ribiere method (5). This technique enhances the steepest descent approach by developing conjugate directions, which ideally minimize the function in fewer iterations. This technique is more appropriate for large-scale problems because of its reduced memory requirements compared to techniques utilizing Hessian matrices. Our study shows its robustness in handling sparse systems, as seen in the highly faster convergence rates compared to the steepest descent in our test cases. However, the performance is widely dependent on the accurate calculation of conjugation conditions, which if not maintained, can lead to suboptimal convergence.

## 3.2. Levenberg-Marquardt Method

The Levenberg-Marquardt method is a numerical method that takes significant steps on steep drops and smaller ones on flat functions. The Levenberg method uses Eq. (6) to consider the function's curvature. The method uses the inverse Hessian Matrix $\nabla^2 U(x^k)$ to read the curve and precisely follow along. The results are precise as well as quick compared to the Steepest Descent.

$$d^{k+1} = -\frac{\nabla U(x^k)}{\nabla^2 U(x^k)} \quad (6)$$

However, like the Steepest Descent method, it slows down as you approach, although less severely than the Steepest Descent. Furthermore, the function can overshoot the minimum point, leaving room for some errors for Levenberg. The steps to the process can be relatively slow since it has to plot the curvature to follow it better. As a result, the Levenberg-Marquardt method takes a long time to estimate the minimum and becomes very expensive. This method finds exact solutions quickly when the function is well-behaved and the initial guess is close to the actual minimum. Our findings show that despite its high computational cost because of Hessian matrix calculations, the method substantially decreases the number of iterations needed to achieve convergence in nonlinear optimization problems. However, its efficiency reduces with poor scaling and non-convex functions, where Hessian calculations can become intractable.

## 3.2. Newton's Method

Newton's method, aka the Newton-Raphson method, is a numerical method with slow steps similar to the Steepest Descent however remains more effective. Newton's method uses the Eq. (7) to take the curvature and slope into consideration. The method uses the Hessian Matrix H, as calculated in (8), a second-order matrix. The results are precise as well as quick compared to the Steepest Descent.

$$d^{k+1} = -H\nabla U(x^k) \tag{7}$$

$$H = [\nabla^2 U(x^k)]^{-1} \tag{8}$$

However, the steps remain huge and slow due to all the calculations that go into this, resulting in this method being very expensive. Not only does it have to plot out itself accurately, but it also must plot the Hessian Matrix in order to do so, which further increases the cost. The resulting program is costly and takes a long time to run, making it very inefficient to run extensive functions to search for the targeted point. Our experiments suggest that the Broyden-Fletcher-Goldfarb-Shanno algorithm is highly effective in medium to large-scale optimization problems, which offers a good compromise between convergence speed and iteration cost. Nevertheless, the success of BFGS hinges on a well-chosen initial approximation and can be sensitive to the accuracy of gradient estimates.

## 3.3. Quasi-Newton's Method

To make the Newton-Raphson method cheaper and faster, the Quasi-Newton method took a slightly different approach. To do this, they approximate the Hessian Matrix using only first-order derivatives, so while it is slower on conversions, it is overall faster than the original. Eq. (7) is still used to calculate the direction; however, depending on the value of k, the Hessian Matrix is set equal to either (9) or (10), where I is the identity matrix starting from the Steepest Descent method.

$$\begin{aligned} &{}^{-1} \text{ when k > 0} \\ &H^k = I \; when \; k = 0 \end{aligned} \tag{9}$$

$$H^{k+1} = H^k + \frac{y^k (y^k)^T}{(y^k)^T s^k} - \frac{[H]^k s^k (s^k)^T [H]^k}{(s^k)^T [H]^k s^k} \tag{10}$$

$$H^{k+1} = H^k + \frac{y^k (y^k)^T}{(y^k)^T s^k} - \frac{[H]^k y^k (y^k)^T [H]^k}{(s^k)^T [H]^k s^k} \tag{11}$$

The first quasi method successfully developed was the DFP, which stands for the Davidson Fletcher Powell method. He had the idea to estimate the Hessian based on the gradient instead of calculating it, which lowered the cost. The approximated Hessian had to be symmetric, the gradient of the model had to equal the function's gradient at $x_k$ and $x_{k-1}$, displacement vector $s_{k-q}= x_k-x_{k-1}$, and gradient change $y_{k-1}=\nabla U(x_k)- \nabla U(x_{k-1})$, furthermore, the resulting matrixes should be as close as possible as shown in Eq. (10). To further optimize this method, the Brayden-Fletcher-Goldfarb-Shannon method Eq. (11) was developed to speed up the process. Rather than approximating the current point, they simplified the process by analyzing the inverse Hessian Matrix, which was faster. It remains the fastest gradient-based method and is less sensitive to step size choice than the others [41]. This technique is specifically valuable in applications like curve fitting and machine learning where the error landscape is typically non-linear. Our analysis shows that it provides stable convergence properties and is robust against poor initializations and ill-conditioned problems. However, the method requires tuning of its damping parameter to balance the trade-off between the descent direction and the step size, which can be challenging to optimize in practice

## 4. Results And Discussion

### 4.1. Problem Statement

In this paper, we discuss classical gradient-based optimization methods. For testing and comparison of the gradient-based methods, we used four famous equations as follows:

$$f_{Rosenbrock}(x,y) = (a-x)^2 + b(y-x^2)^2 \quad (12)$$

In mathematical optimization, the Rosenbrock function is a non-convex function used as a performance test problem for optimization algorithms. The global minimum is inside a long, narrow, parabolic-shaped flat valley. To find the valley is trivial. To converge to the global minimum, however, is difficult. For this testing problem, $a = 1$ and $b = 100$ are considered.

$$f_{SpringForce}(x,y) = \frac{1}{2}k_1\left(\sqrt{x^2+(l_1-y)^2} - l_1\right)^2 + \frac{1}{2}k_2\left(\sqrt{x^2+(l_2-y)^2} - l_2\right)^2 - p_1x - p_2y \quad (13)$$

Vander Plaat's unconstrained minimization problem is a classic two-dimensional optimization problem. One response function has been chosen chiefly to illustrate how the metamodel assembly works. The objective is to find an equilibrium position of the springs by minimizing the total potential energy of the system. The constants $K_i$ are spring stiffnesses, $P_i$ are loads, $l_i$ are the original spring lengths, and $x_i$ are displacements. For this problem $K_1$ = 8 N/cm, $K_2 = 1$ N/cm, $P_1 = P_2 = 5$ N, $l_1 = l_2 = 10$ cm.

$$f_{Ackley}(x,y) = -20e^{-0.02\sqrt{\frac{x^2+y^2}{2}}} - e^{\frac{cos(2\pi x)+cos(2\pi y)}{2}} + e + 20 \quad (14)$$

The Ackley function is a non-convex function used as a performance test problem for optimization algorithms in mathematical optimization. The Ackley function is widely used for testing optimization algorithms. Its two-dimensional form is characterized by a nearly flat outer region and a large hole at the center. The function poses a risk for optimization algorithms, particularly hill-climbing algorithms, to be trapped in one of its many local minima.

$$f_{Himmelblau}(x,y) = f = (x^2 + y - 11)^2 + (x + y^2 - 7)^2 \quad (15)$$

In mathematical optimization, Himmelblau's function is a multi-modal function used to test the performance of optimization algorithms. It has one local maximum at and where and four identical local minima. The locations of all the minima can be found analytically. However, because they are roots of cubic polynomials, the expressions are somewhat complicated when written in terms of radicals. The following minimization functions are tested using the classical gradient-based methods. The comparison is made with the computing efficiency and accuracy of each of these methods.

## 4.2. Results of Optimization Methods

In this paper, we used six gradient-based methods to optimize four famous minimization problems. The presented methods are the steepest descent, conjugate gradient methods of Fletcher-Reeves and Polak-Ribiere, Newton-Raphson, Quasi-newton (BFGS), and Levenberg-Marquardt methods. The optimization is done for Rosenbrock, Spring force Vanderplaats, Ackley, and Himmelblau's functions, which have been presented in the method section. Results of the minimization problem to find the nearest local minimum are illustrated in the following Figures 1 and 2.

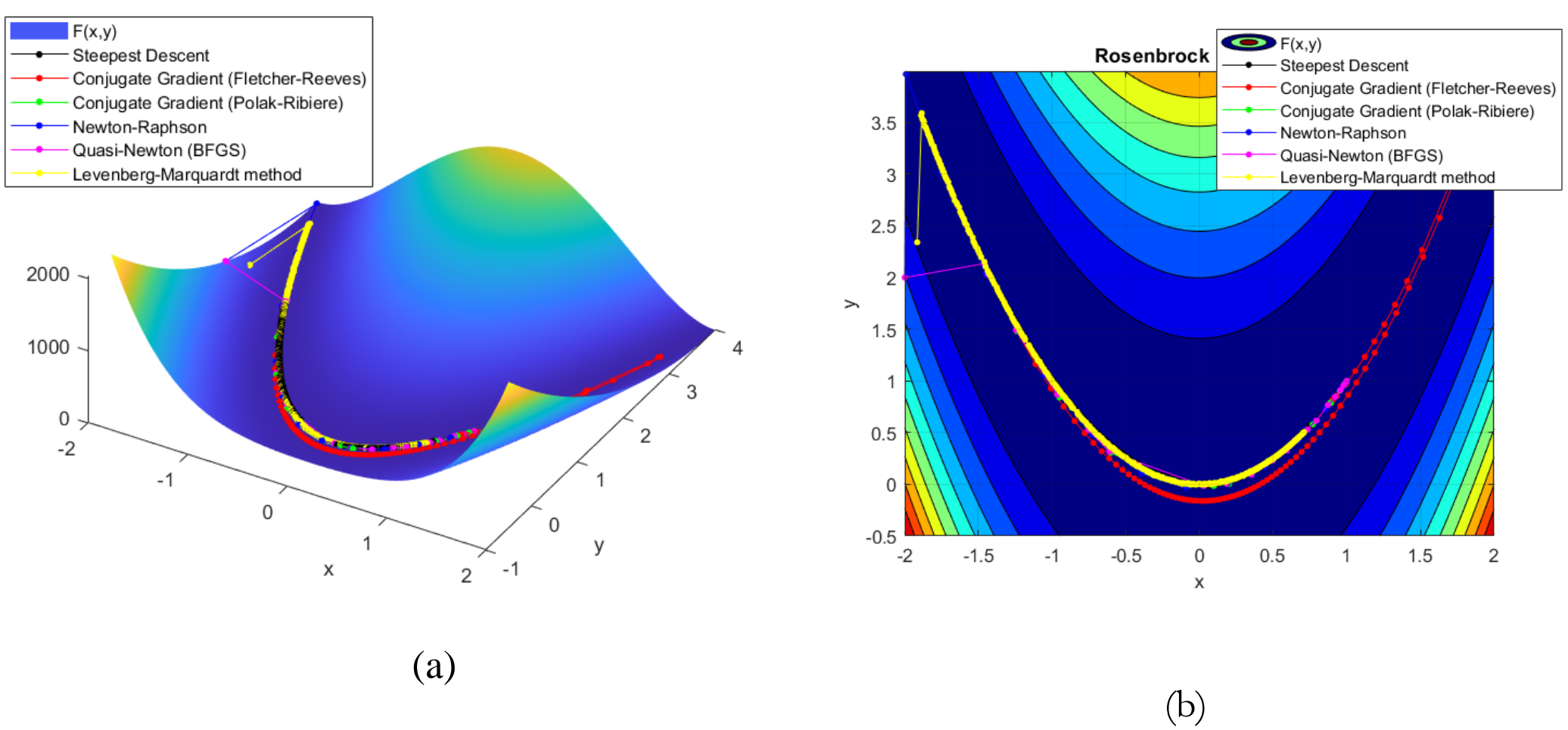

(a)

(b)

**Figure 1:** The results of optimization of Rosenbrock function (a) 3D surface, (b), Contour plot

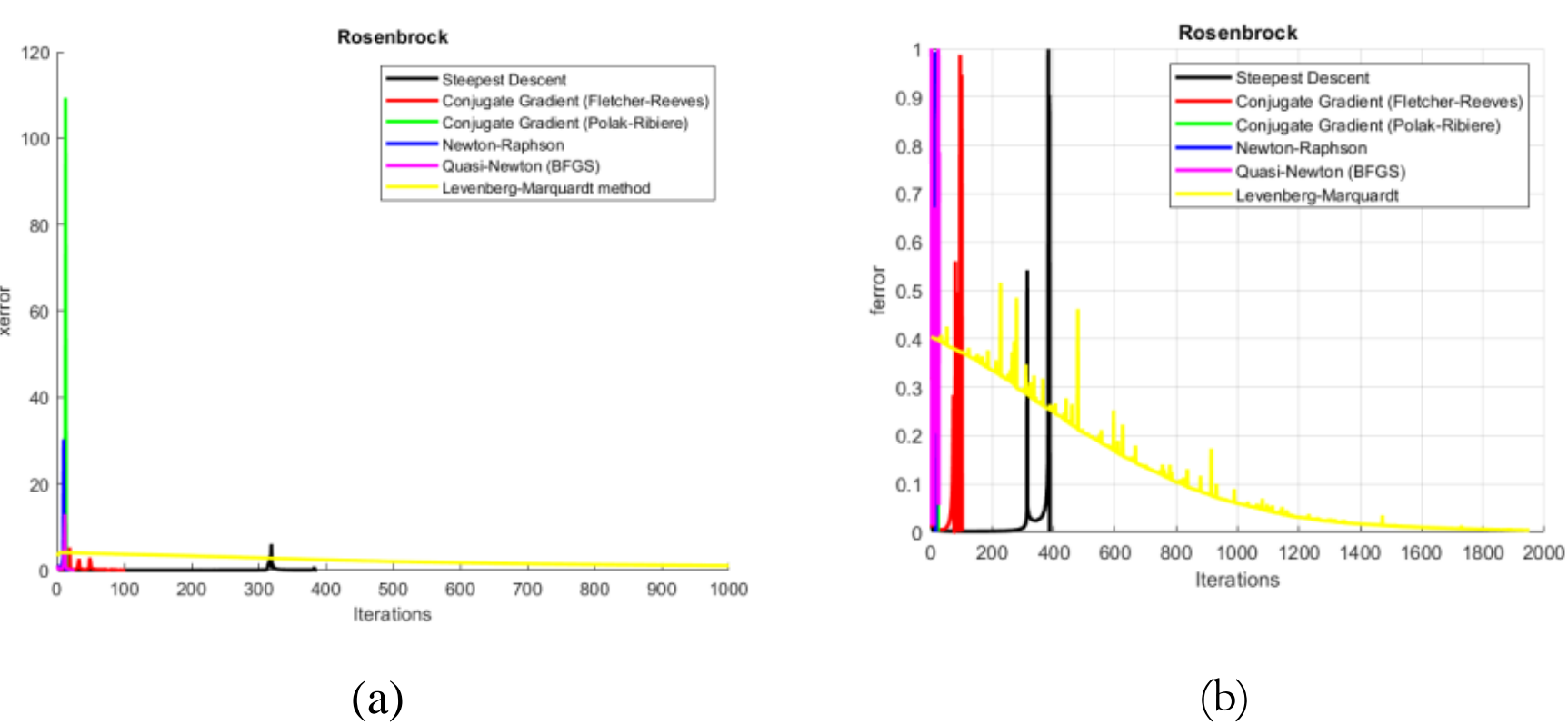

(a)

(b)

**Figure 2**: The results of optimization of Rosenbrock function, (a) the Error value of distance to evaluate the position and final optimum point, (b) the value of function error based on the difference between the value of each iteration and final objective function value

Figure 1 shows the surface and contour plot of the Rosenbrock function. The optimization environment includes $x \in [-2,2]$ and $y \in [-1,4]$. The location of the local minimum for this interval is $(x_m, y_m) = (1,1)$. Moreover, the initial search point is $(x_0, y_0) = (-2,2)$. Based on the results of the optimization, the minimum point with a small error is $(x_p, y_p) = (0.9992, 0.9984)$. Regarding the error of the evaluation process for the x value that is shown in Figure 2 (a), it can be estimated that the greatest iteration or slowest method belongs to the Levenberg-Marquart method with 2000 iterations to stop with $10^{-6}$ stopping criteria. Moreover, the second slow method is the steepest descent approach, with almost 400 process rounds. The error of the Rosenbrock function value is depicted in Figure 2(b). Based on the optimization process results, the presented algorithms see some fluctuation of tolerance in the optimum point that comes back to find the minimum. Regarding the finding of this paper based on Figure 2, the fastest method for the Rosenblock function problem is Newton-Raphson and BFGS algorithms.

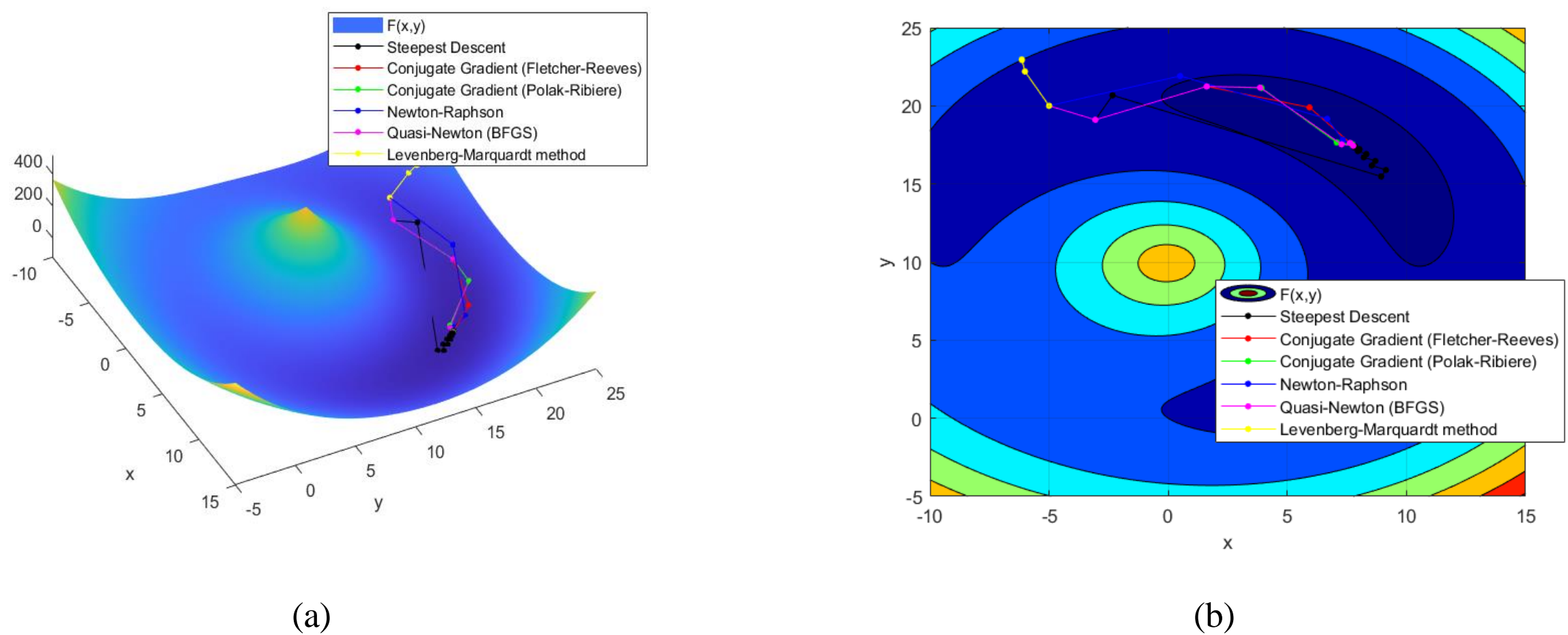

(a) (b)

**Figure 3**. The results of optimization of Spring force function (a) 3D surface, (b), Contour plot

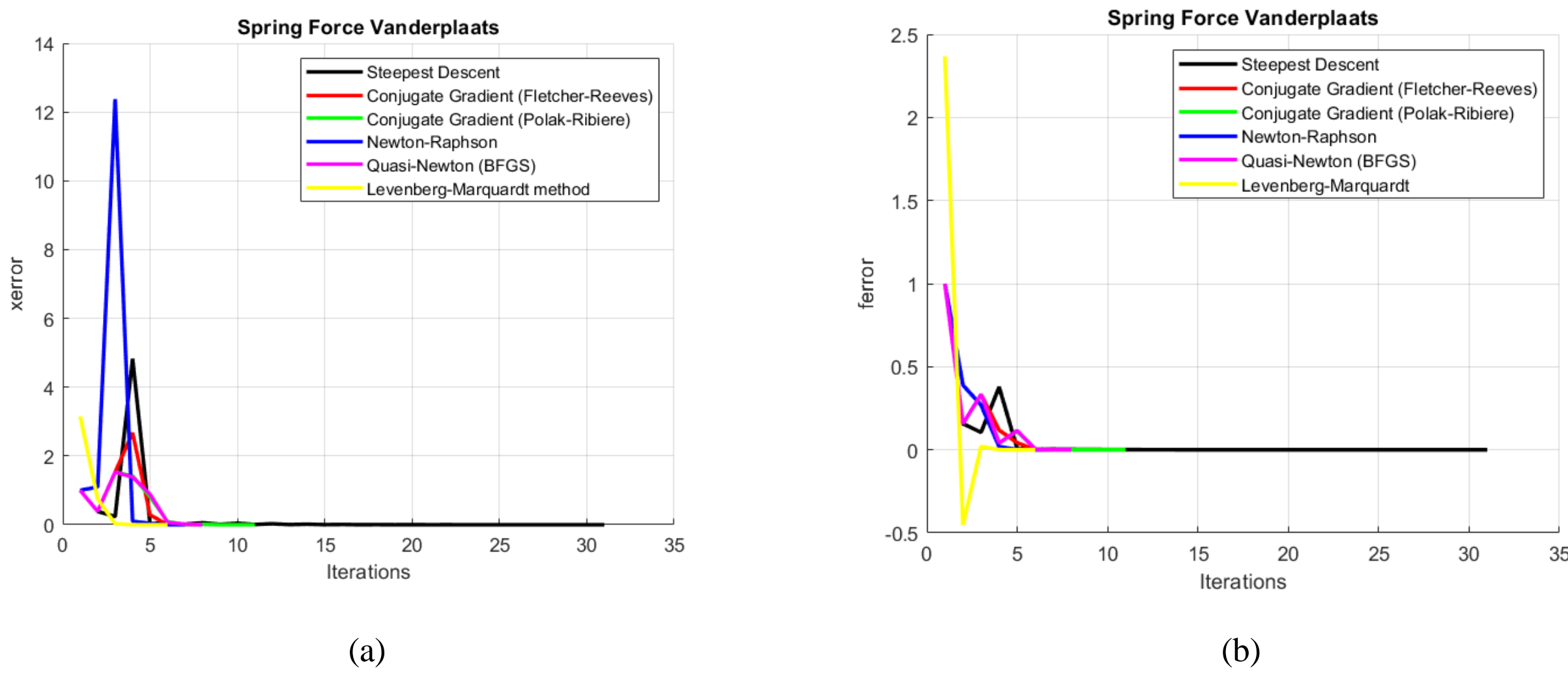

(a) (b)

**Figure 4.** The results of optimization of the Spring force function, (a) the Error value of distance to evaluate the position and final optimum point, (b) the value of function error based on the difference between the value of each iteration and final objective function value

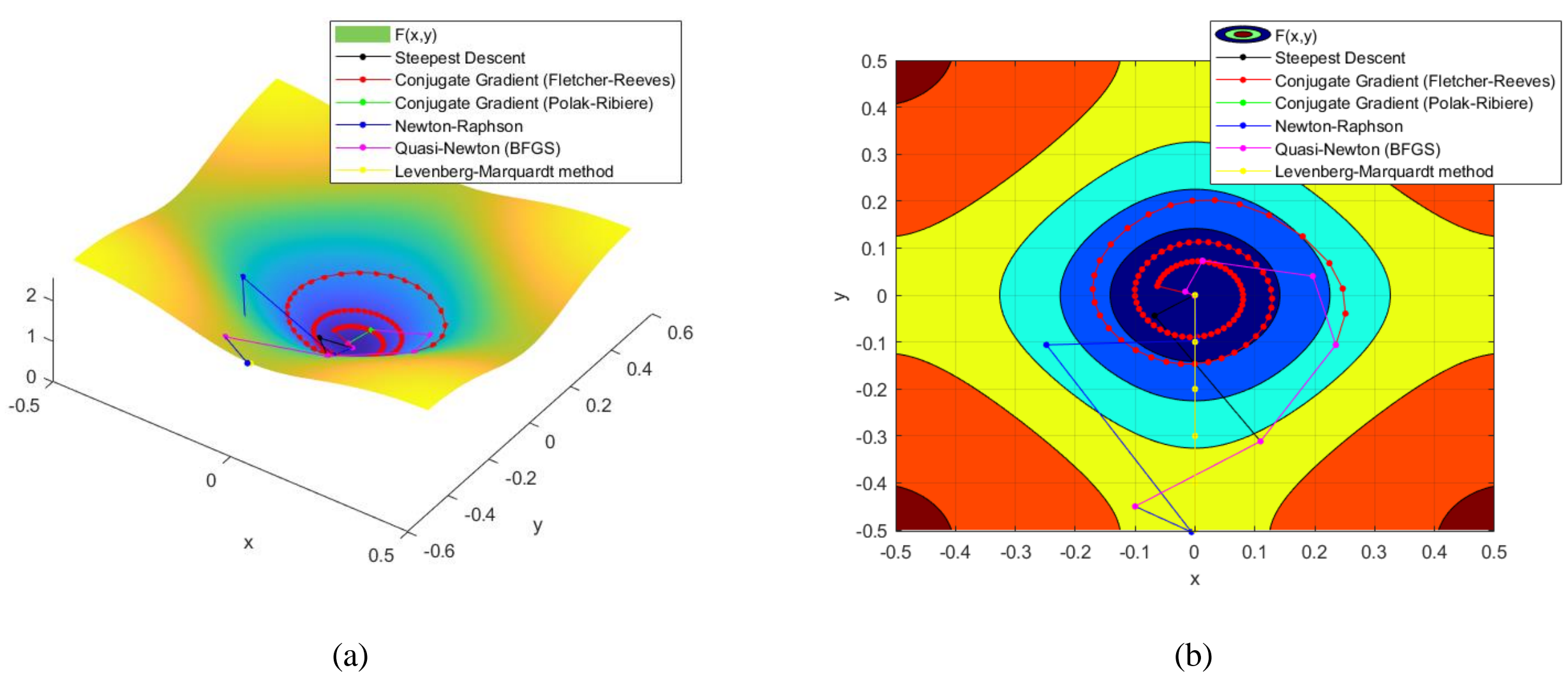

**Figure 5**. The results of optimization of Ackley function (a) 3D surface, (b), Contour plot

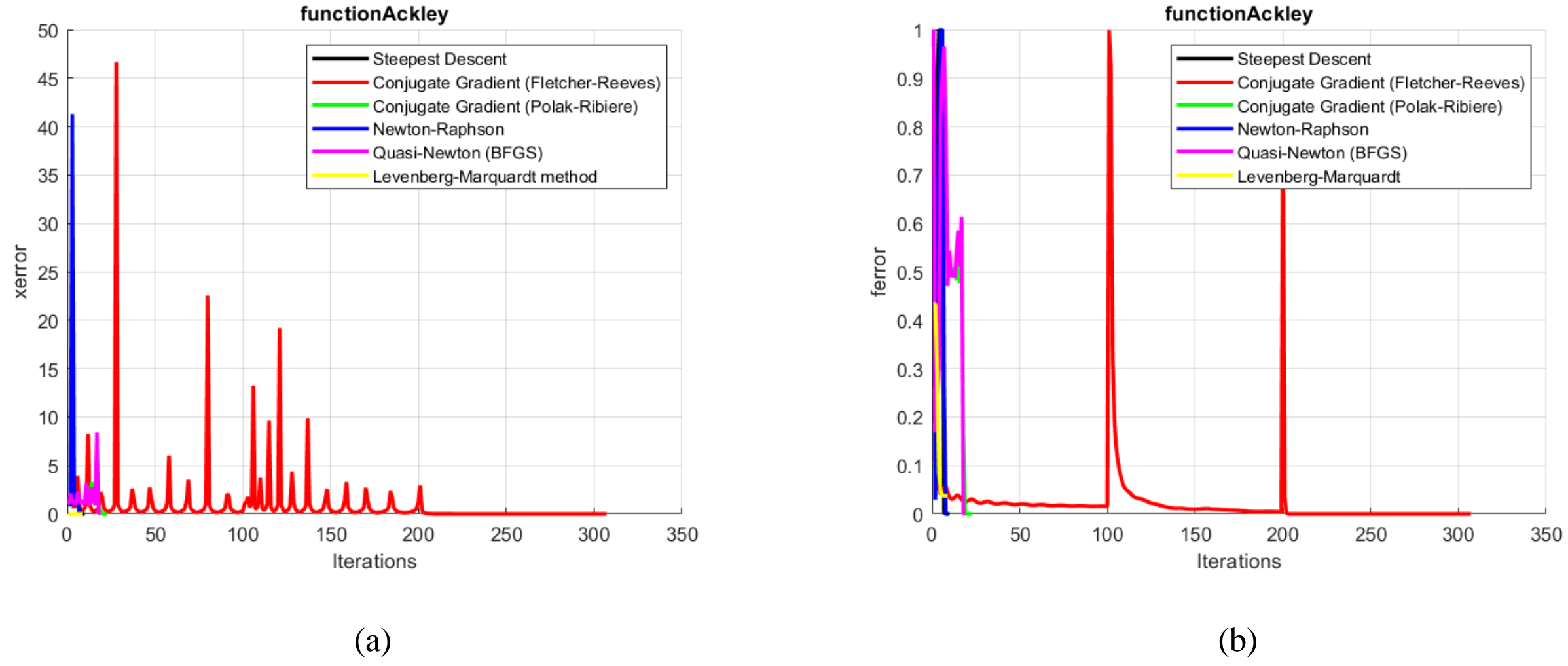

**Figure 6:** The results of optimization of Ackley function, (a) the Error value of distance to evaluate the position and final optimum point, (b) the value of function error based on the difference between the value of each iteration and final objective function value

The surface and contour plots of the Spring Force function are shown in Figure 3. The optimization environment includes x∈[-10,15] and y∈[-5,25]. The location of the local minimum for this interval is (x_m,y_m )=(7.623,17.626). Moreover, the initial search point is (x_0,y_0 )=(-1,1). The minimal point with a modest error is (x_p,y_p )=(7.626930,17.626331) based on the optimization findings. By calculating the error of the assessment process for the x value shown in Figure 4 (a), the steepest descent method has the most iterations or is the slowest method, with 31 iterations to stop with and 10-6 terminating criteria. Moreover, the conjugate gradient Polak-Riviere method, which has about ten rounds in the process, is the second slow method that can be applied. Figure 4(b) shows the error of the Spring Force function value.

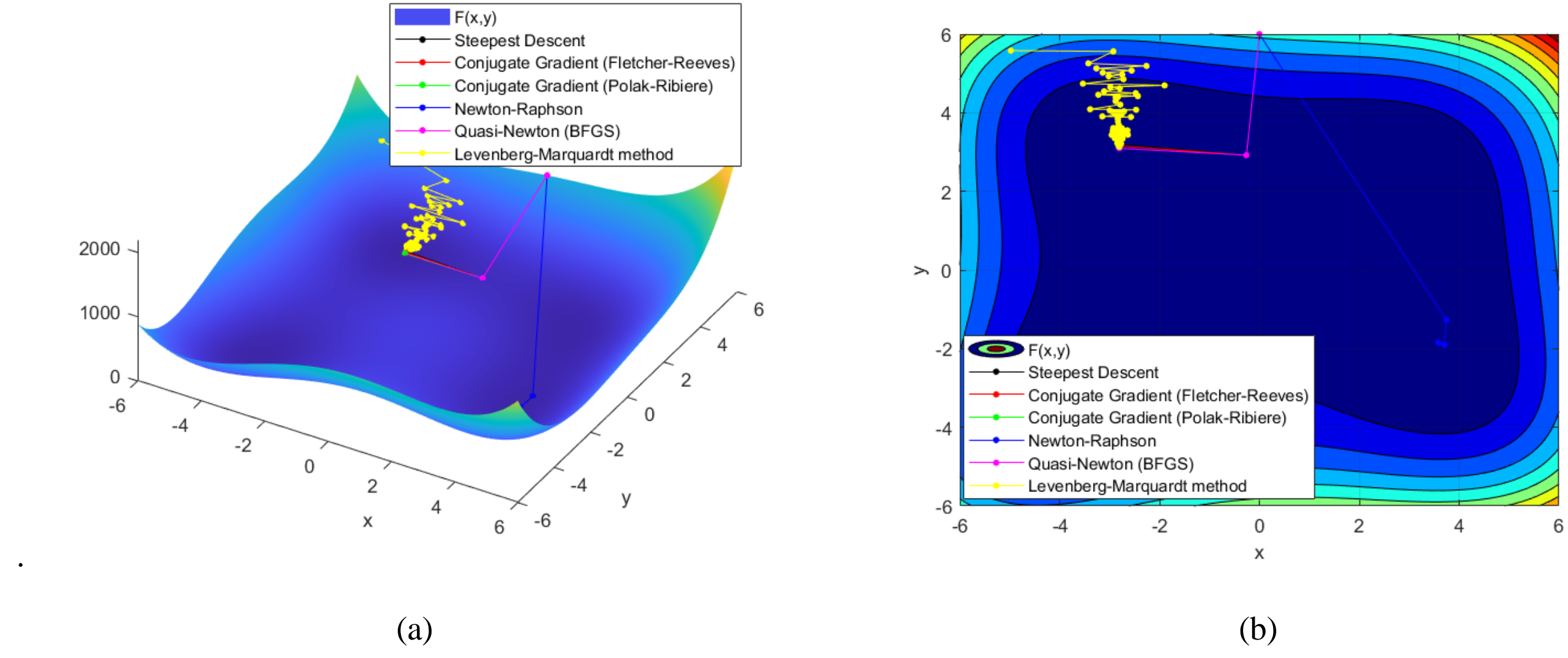

(a) (b)

**Figure 7**. The results of optimization of Himmelblau's function (a) 3D surface, (b), Contour plot

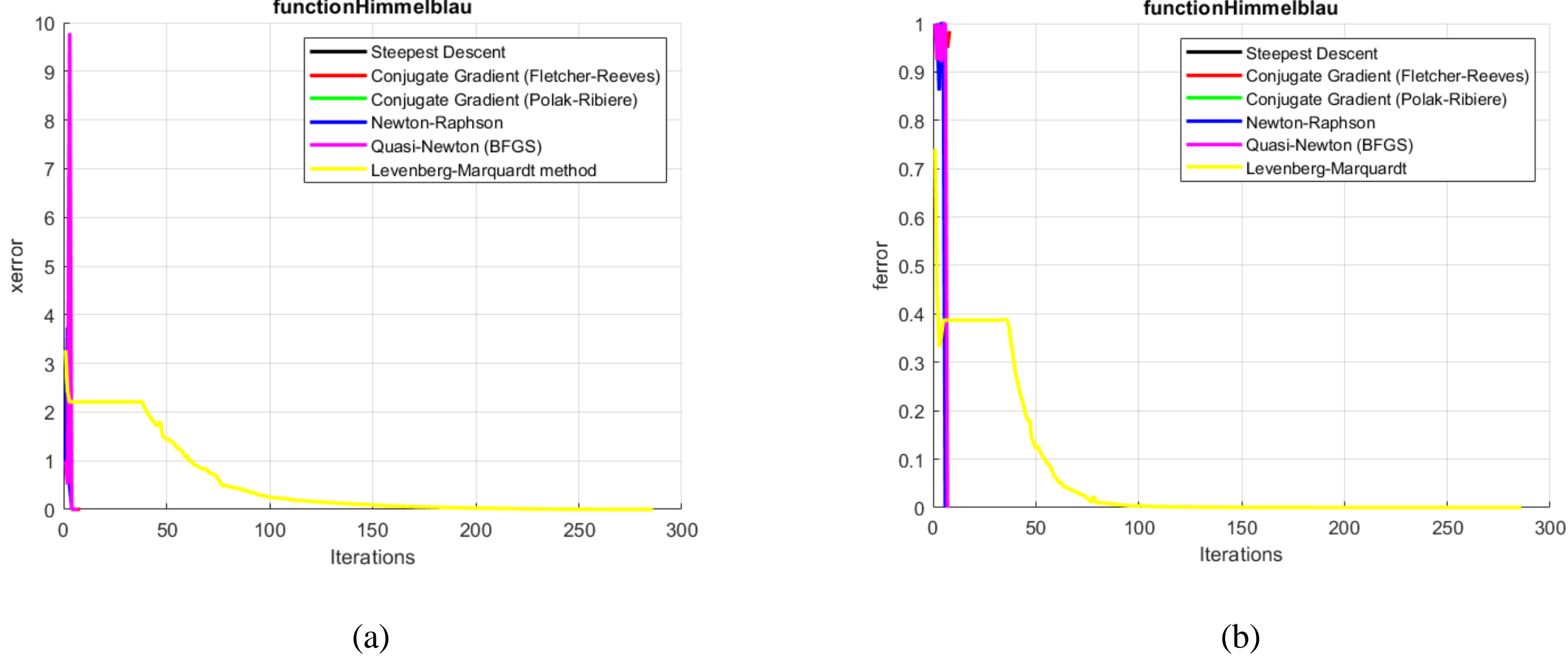

(a) (b)

**Figure 8**: The results of optimization of Himmelblau's function, (a) the Error value of distance to evaluate the position and final optimum point, (b) the value of function error based on the difference between the value of each iteration and final objective function value

According to the findings of this study based on Fig. 4, Newton-Raphson is the quickest way to solve the Spring Force function problem. The Levenberg-Marquart algorithm cannot locate the minimal points with the most fantastic accuracy while optimizing this equation. Figure 5 shows a plot of the Ackley function and a contour plot of the function. The calculation interval is $x\in[-0.5,0.5]$, and $y\in[-0.5,0.5]$. The position of the local minimum for this rectangular interval is $(x_m,y_m)=(0,0)$. Moreover, the initial search point is $(x_0,y_0)=(-0.1,-0.45)$. The presented gradient-based methods can find the extremum point with the highest accuracy. Based on the error of the assessment process displayed in Figure 6 (a), we can conclude that the conjugate gradient Fletcher-Reeves technique with 310 iterations is the slowest technique. The inaccuracy of the Ackley function value can be seen in Figure 6(b). Using Figure 6, the study suggests Levenberg-Marquart be the fastest method for solving the Ackley function problem.

Figure 7 shows a plot of Himmelblau's function and a contour plot of the function. The calculation interval is $x\in[-6,6]$, and $y\in[-6,6]$. This function has four local minimums in this interval. The initial search point is (x_0,y_0 )=(0,6). Regarding Figure 7, Newton-Raphson's methods cannot find the nearest local minimum and show another local minimum that is far from the initial point. However, other methods find the nearest local extremum with high accuracy. We conclude that the Levenberg-Marquart technique, with 280 iterations, is the slowest among all alternative methods based on the error rate shown in Figure 8 (a).

## 4.3. Effects of Initial Point Position of Optimization

We study the effect of initial points position on optimization methods. For this purpose, 100 initial points at a rectangle with a 10x10 size near the optimum point are selected. For the study of the conjugate gradient method of Fletcher-Reeves, the Rosen brock function is used. The results of optimization are presented in Figure 9. Based on the results, the selected initial points have no significant effect on the optimum point. The algorithm can optimize the function with almost the highest accuracy. Moreover, finding the minimum point is almost similar in all cases. This study's results for the steepest descend method are also shown in Figure 10. Like the process is done for all conjugate gradients, the Polak-Ribiere method with spring force function (Figure 11), Newton-Raphson with Himmelblau's function (Figure 12), BFGS with Ackley function (Figure 13), and Levenberg-Marquardt with Rosenbrock function (Figure 14). Regarding the results, the effect of initial points for steepest descends, conjugate gradients Fletcher-Reeves, and conjugate gradients Polak-Ribiere Levenberg-Marquardt and DFGS are insignificant. In the Newton-Raphson method, some of the points found cannot find the optimum because of the complexity of Himmelblau's with triple minimum points. Regarding Fig. 12, there is no unique x and y value to find optimum points.

The values are concentrated on a curve instead of the optimum point. Building upon these findings, this study expands the discussion to encapsulate the wider importance and theoretical improvements of gradient-based optimization techniques in various computational settings. Our analysis emphasizes the adaptability of these techniques in addressing non-linear and complex optimization problems, which offers a transparent picture of their performance across multiple scenarios described by varying degrees of multimodality and constraint intricacies. The results derived from this work contribute to the field by suggesting an integrated practical framework that emphasizes the potential of combining traditional gradient-based techniques with state-of-the-art optimization methods. Our approach advances discourse within the optimization scope. It also positions the stage for more advanced and transformative applications in machine learning as well as engineering sectors, where solving dynamic and multi-dimensional problems is vital.

By investigating the sensitivity to initial conditions, comparing diverse methods, and analyzing their effectiveness across multi-modal and non-linear functions, this study presents important viewpoints for improving the selection of optimization algorithms. Also, the discussion on incorporating heuristic methods opens ways for tackling complicated optimization tasks that traditional approaches struggle to address alone. These contributions improve the current frameworks and offer practical advantages to industries that rely on advanced optimization, from engineering to tech-driven sectors, which shows the importance of the development of adaptable and more robust optimization solutions.

In addition to the influence of initial points on the optimization performance, it is vital to consider the roles of other essential hyperparameters, including the step size (often referred to as the learning rate in machine learning contexts) and regularization factors. The step size $\alpha$ is a critical control parameter in the convergence process of gradient-based techniques. An appropriately chosen step size ensures that the optimization process is neither too slow (small step size leading to slow convergence) nor too unstable (large step size causing overshooting of the minimum). Various adaptive step size techniques can be employed to dynamically adjust the learning rate based on the iteration progress, such as learning rate schedules or adaptive gradient methods like AdaGrad, RMSprop, and Adam. Moreover, regularization plays a notable role, particularly in scenarios prone to overfitting or where the objective function is ill-posed. Regularization techniques, such as L1 and L2 regularization, add a penalty term to the objective function. This penalty helps control the model complexity, ensuring that the solution does not fit the noise in the data. The choice of regularization factor $\lambda$ can significantly affect the optimization's outcome by balancing the trade-off between bias and variance in the model.

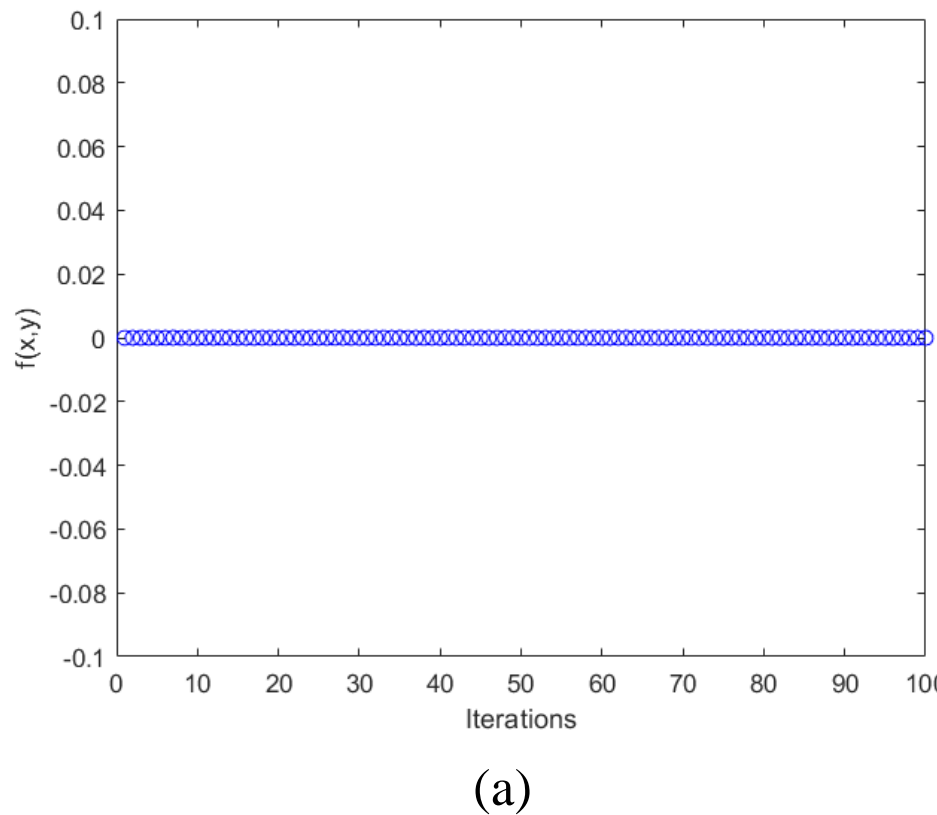

(a)

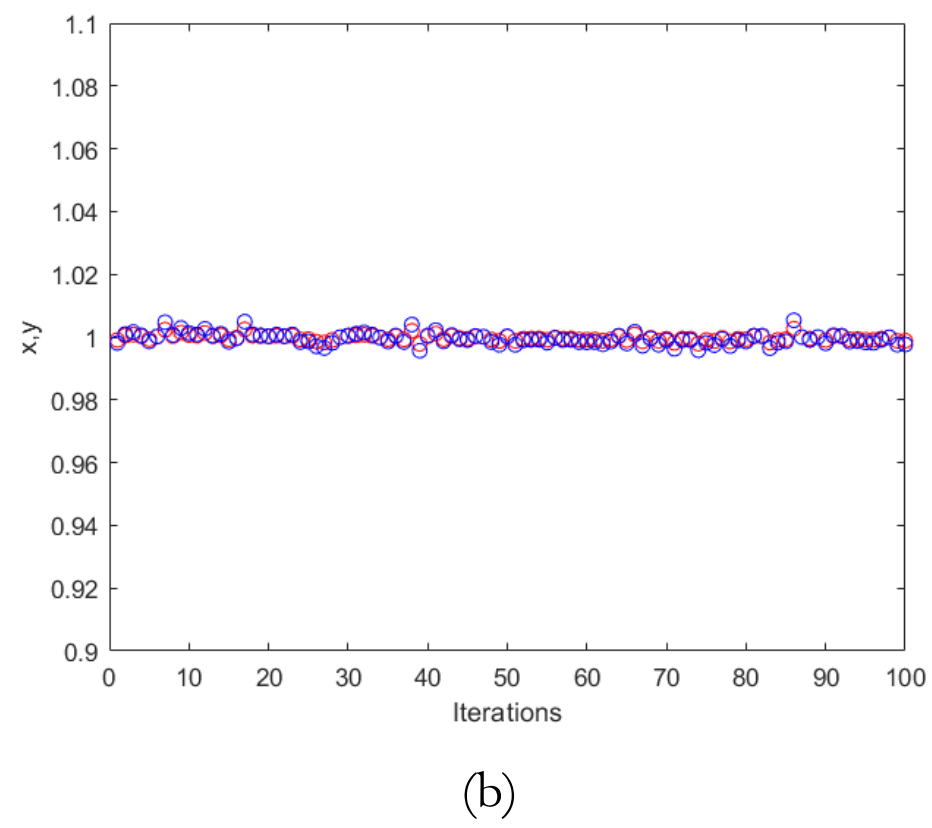

(b)

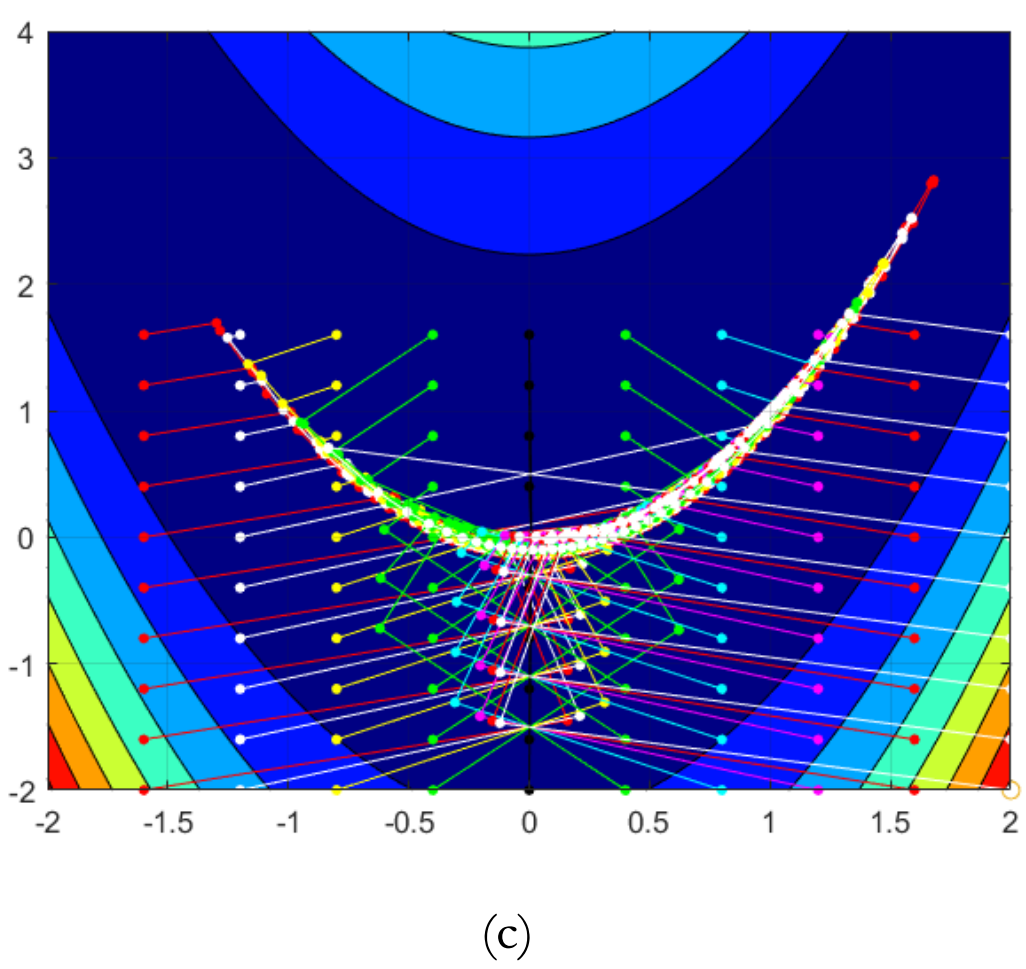

(c)

**Figure 9**. Different initial points in optimization for conjugate gradients Fletcher-Reeves with Rosenbrock function(a): optimum function value. (b): The x and y value of optimum points, (c): The process of finding optimum points by starting different initial points

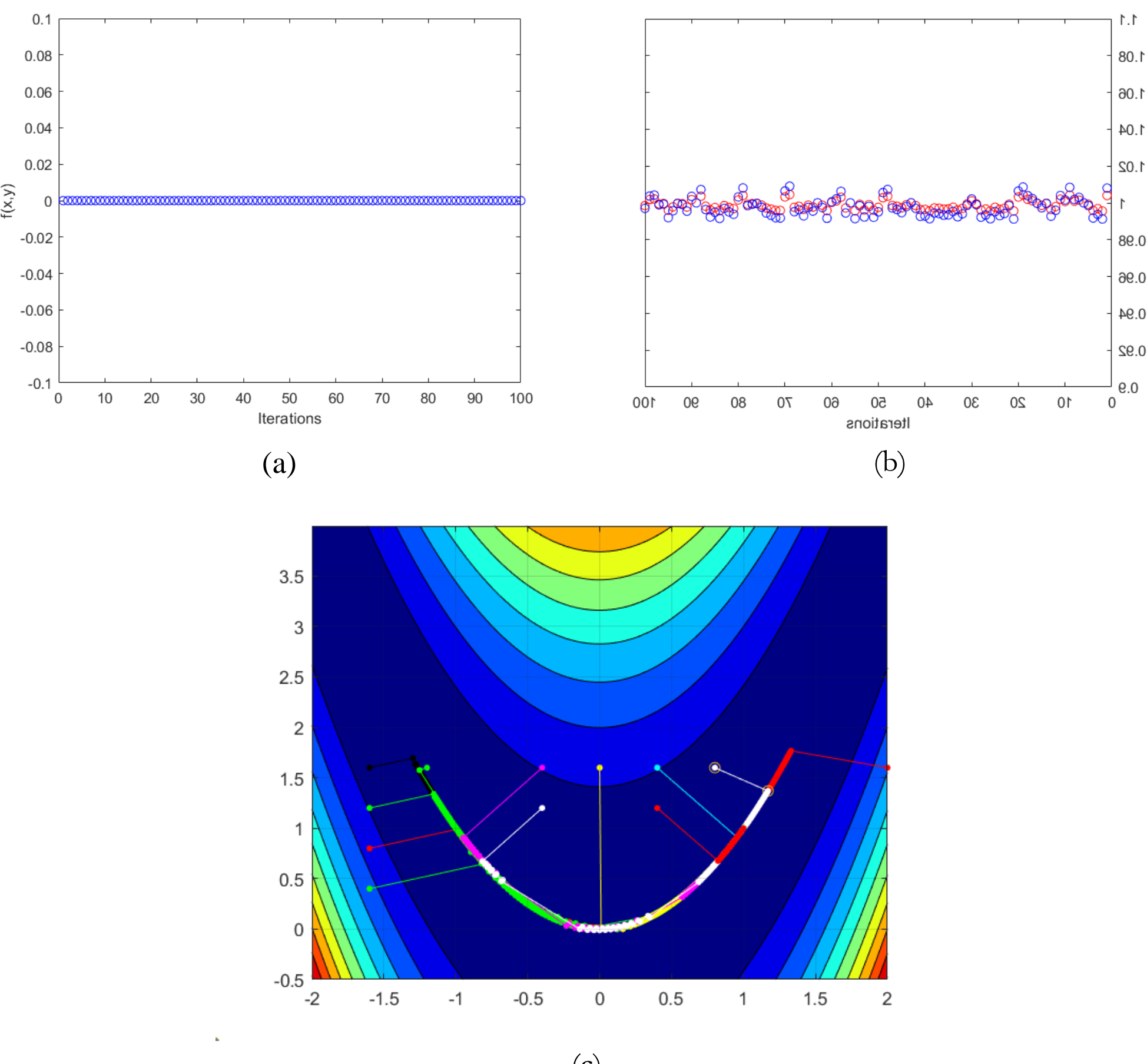

(a)

(b)

(c)

**Figure 10**. Results of different initial points in optimization for Steepest descend with Rosenbrock function(a): optimum function value. (b): The x and y value of optimum points, (c): The process of finding optimum points by starting different initial points

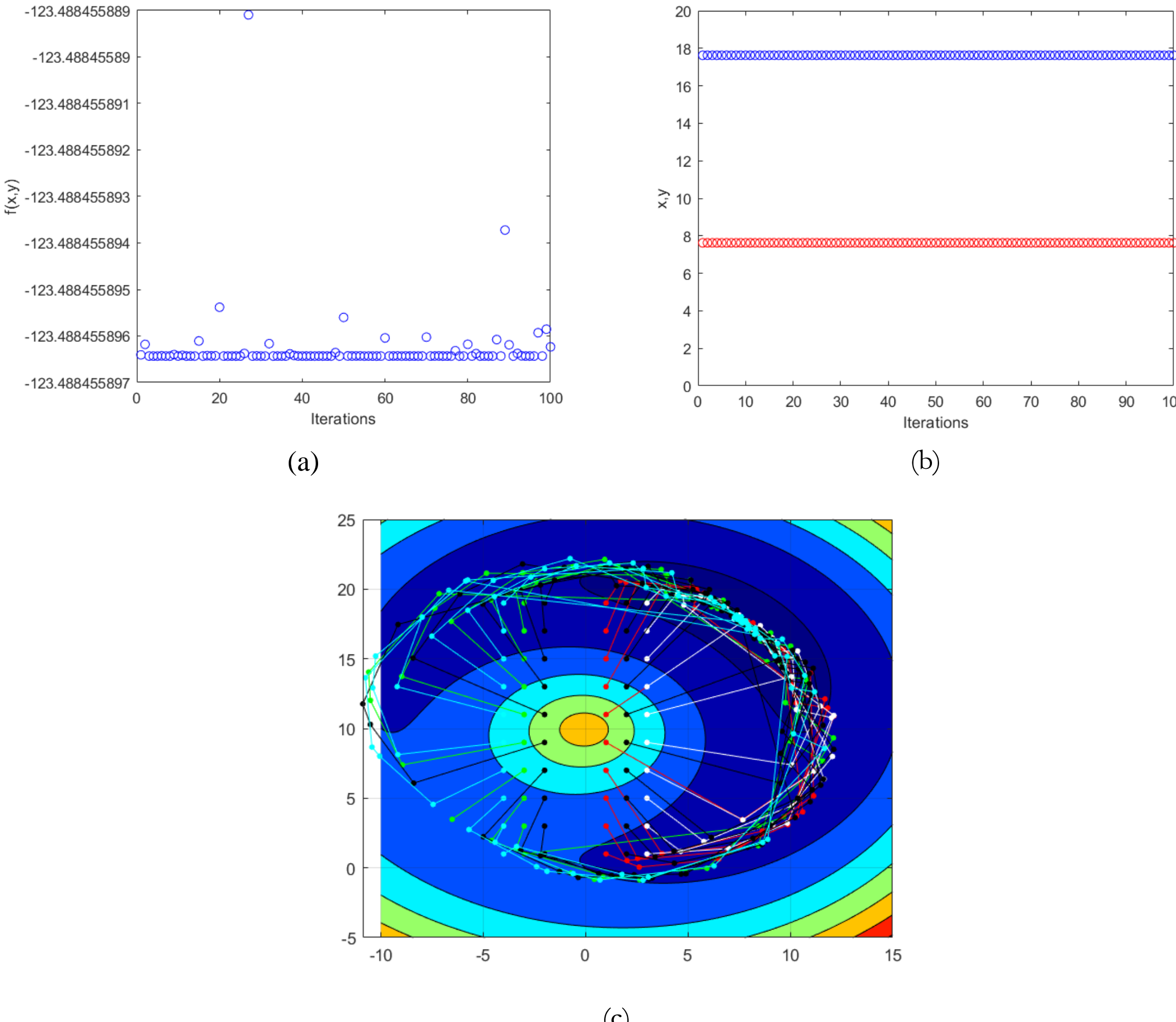

**Figure 11:** Results of different initial points in optimization for conjugate gradients Polak-Ribiere with spring force function(a): optimum function value. (b): The x and y value of optimum points, (c): The process of finding optimum points with starting different initial points

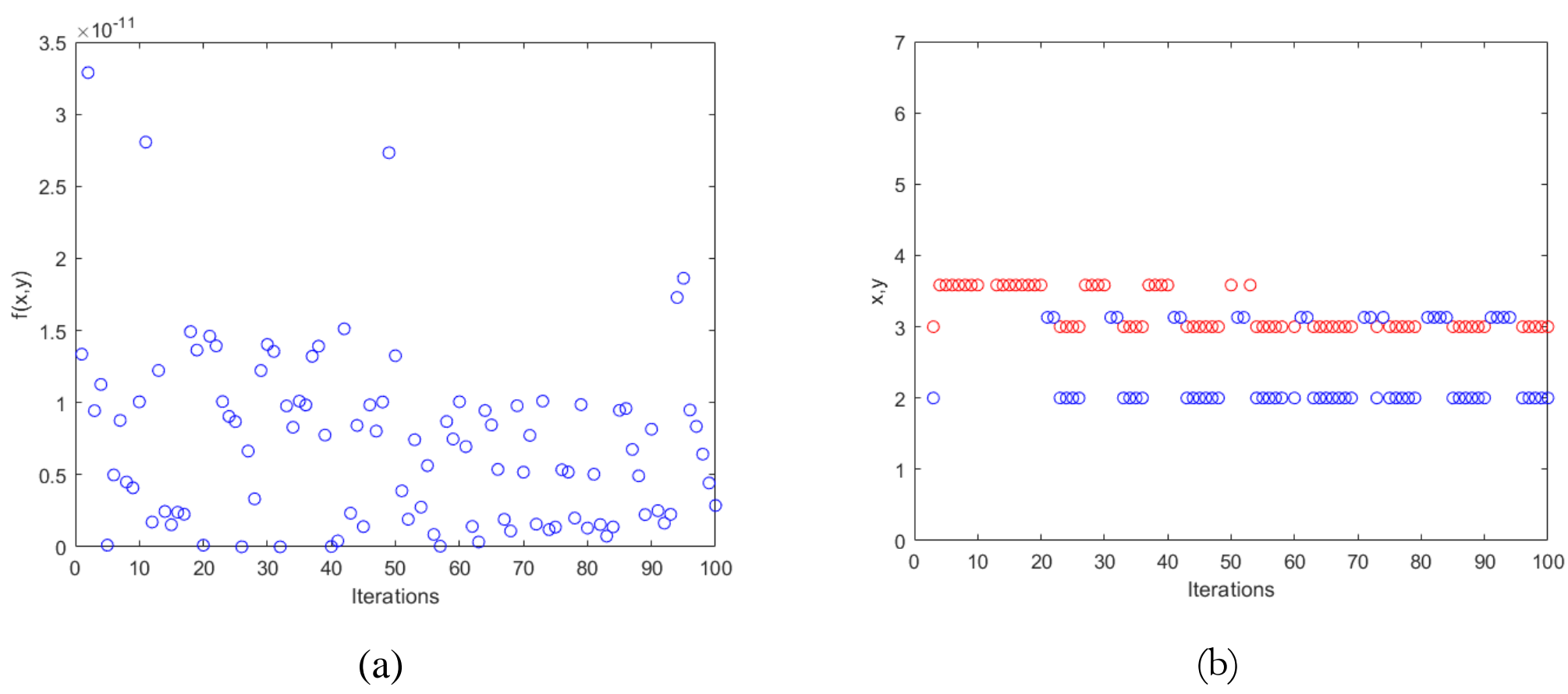

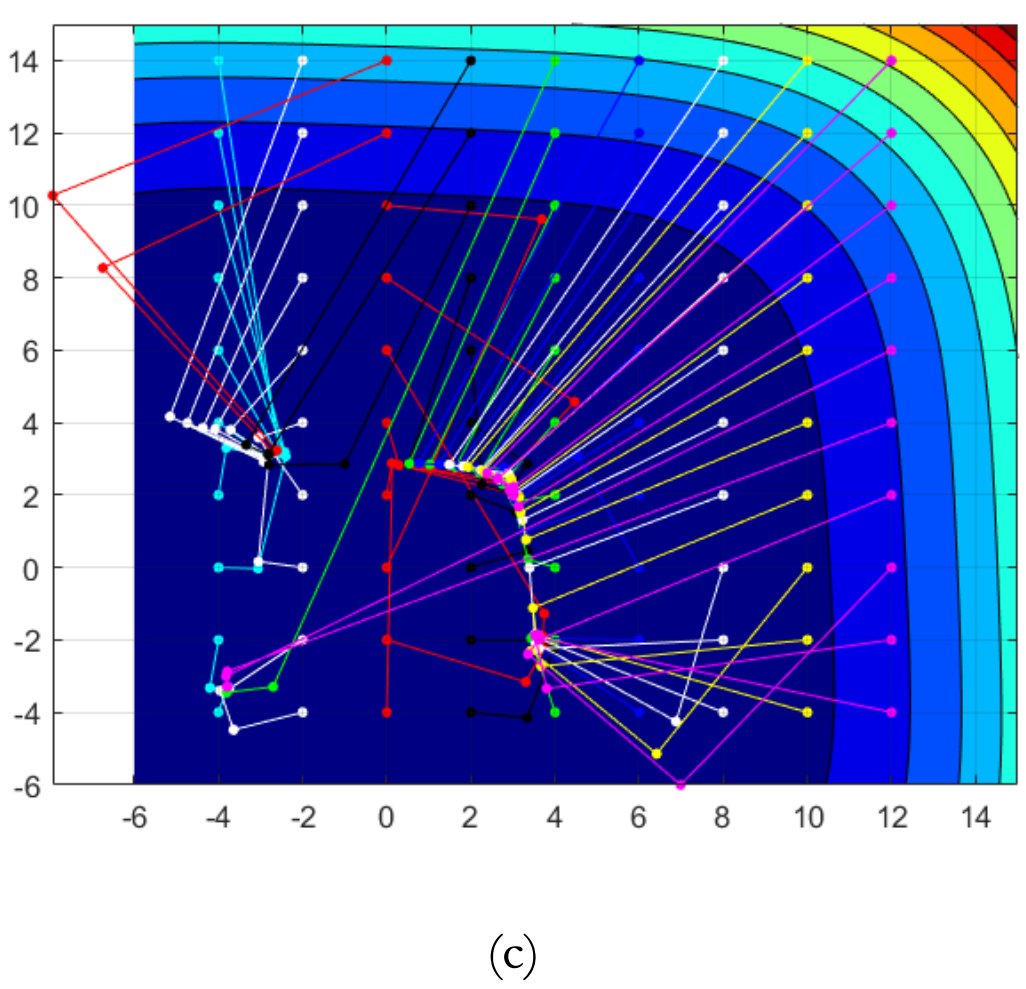

(c)

**Figure 12:** Results of different initial points in optimization for Newton-Raphson with Himmelblau's function(a): The x and y values of optimum points. (b): optimum function value, (c): The process of finding optimum points by starting different initial points

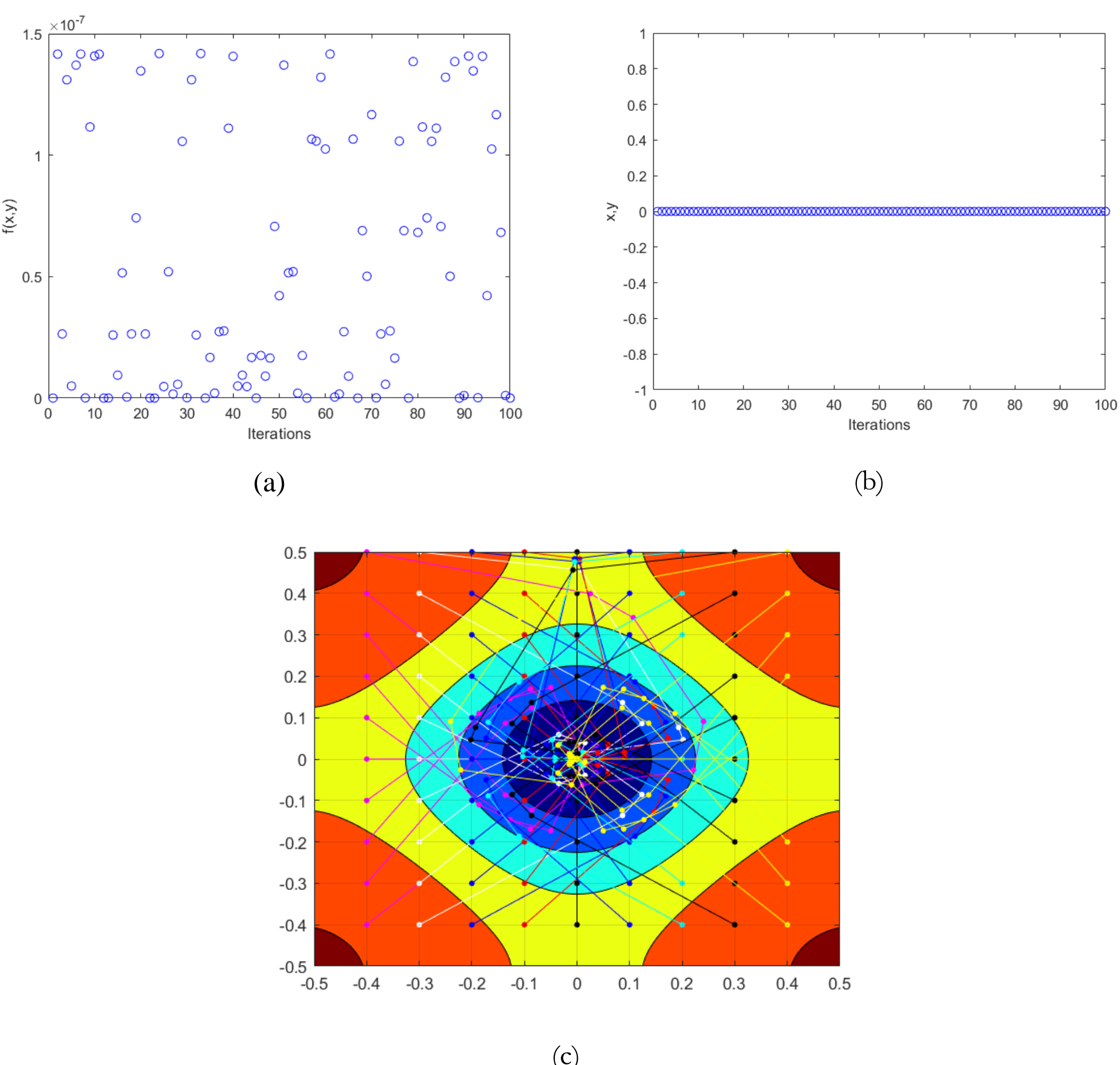

(a) (b)

(c)

**Figure 14**: Results of different initial points in optimization for BFGS with Ackley function(a): The x and y value of optimum points. (b): optimum function value, (c): The process of finding optimum points by starting different initial points

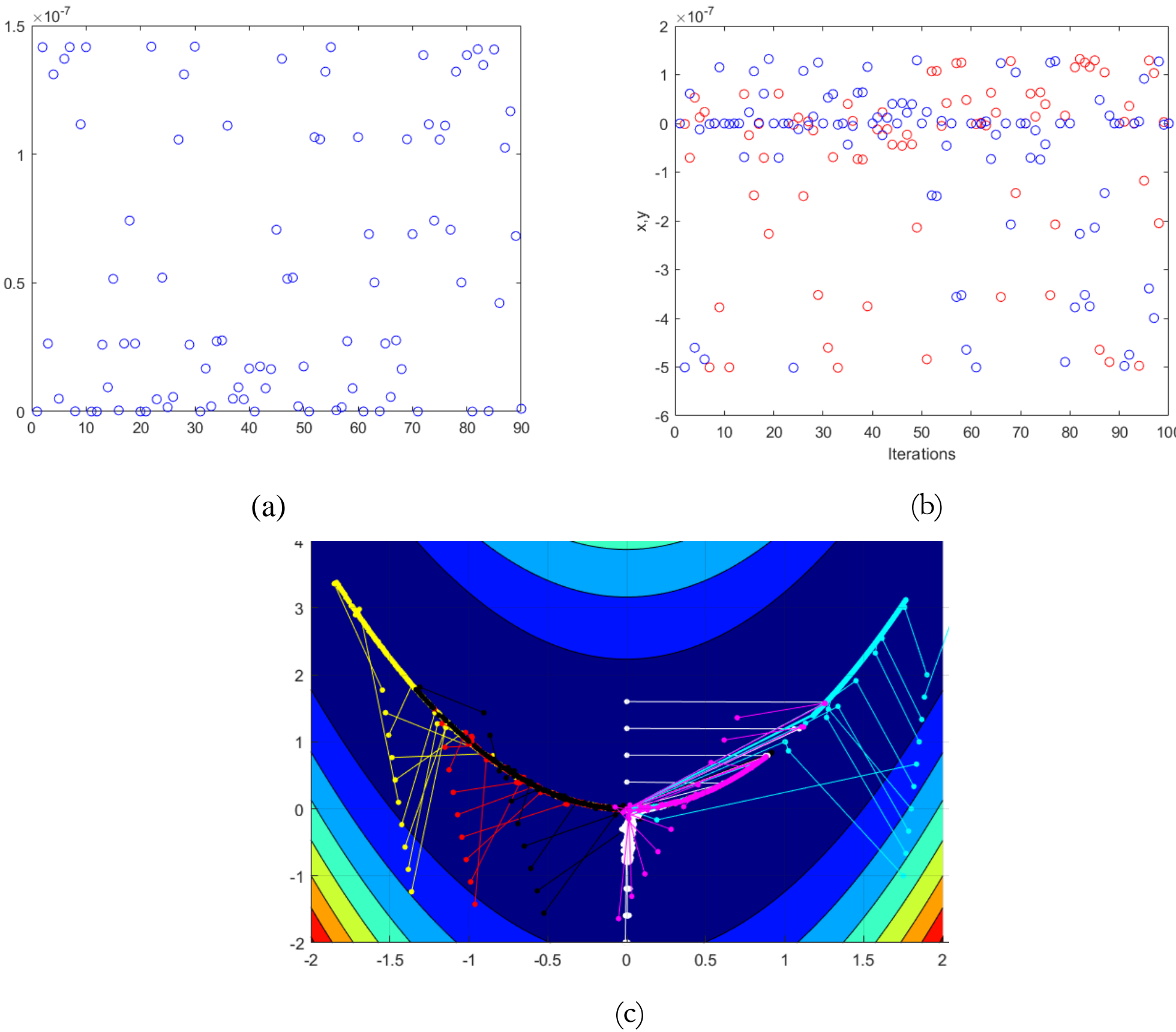

**Figure 15**: Results of different initial points in optimization for Levenberg-Marquardt with Rosenbrock function(a): The x and y value of optimum points. (b): optimum function value, (c): The process of finding optimum points by starting different initial point

## 5. Conclusions

Our project involved several gradient-based optimization algorithms on functions with n-variables and no constraints. In our study of the techniques, we looked at the complexity, accuracy, number of iterations, efficiency, and overall effectiveness of each to determine how efficient it was. We tested five methods on three different functions for the most accurate results. We present the steepest descent, conjugate gradient methods associated with Fletcher-Reeves and Polak-Ribiere, Newton-Raphson, quasi-newton (BFGS), and Levenberg-Marquardt methods. The functions presented in the method section are Rosenbrock, Spring Force Vanderplaats', Ackley's, and Himmelblau's. The choice of the initial guess is also taken into account when evaluating results. Lastly, the optimization process is plotted and evaluated for different 3D surfaces. It is essential to make sure you are looking for the global optimum rather than a local optimum, a global suboptimal solution.

The fundamental reason why traditional approaches find poor answers is because of nonlinearity and multimodality. A gradient-based approach will not operate well if the goal function discontinues. The number of choice factors might also be an issue. The number of possible solutions, along with non-linearities, can surpass most computers' computational capabilities, making searching for all conceivable combinations impossible. This difficulty is addressed by heuristic and metaheuristic algorithms. All optimization problems are written so that the goal and constraint functions may be precisely specified. In reality, all measured variables are subject to some degree of uncertainty. When uncertainty and noise enter the equation, the optimization becomes a stochastic optimization issue, also known as robust optimization with noise. Standard optimization approaches can only be employed if the issue is redefined or recast

**Author Contributions:** All authors have the same contribution

**Funding:** This research received no external funding

**Data Availability Statement:** Data is available and can be provided over the emails querying directly to the author at the corresponding author.

**Conflicts of Interest:** The authors declare no conflicts of interest